\documentclass[epsfig,latexsym,amsfonts,twoside]{article}
\usepackage{amssymb,amsmath,amscd,epsfig}
\usepackage[all,cmtip]{xy}

\def\part#1{\frac{\partial\phantom{#1}}{\partial#1}}
\newtheorem{thm}{Theorem}
\newtheorem{prp}[thm]{Proposition}
\newtheorem{lem}[thm]{Lemma}
\newtheorem{cor}[thm]{Corollary}

\newenvironment{prf}{\begin{trivlist}\item[]{\bf Proof} }%
{\hfill $\Box$ \end{trivlist}}
\newenvironment{dfn}{\begin{trivlist}\item[]{\bf Definition}\em }%
{\end{trivlist}}
\newenvironment{rmk}{\begin{trivlist}\item[]{\bf Remark} }%
{\end{trivlist}}
\newenvironment{exm}{\begin{trivlist}\item[]{\bf Example} }%
{\end{trivlist}}
\newenvironment{que}{\begin{trivlist}\item[]{\bf Question} }%
{\end{trivlist}}

\def\Z{\ifmmode{{\mathbb Z}}\else{${\mathbb Z}$}\fi}
\def\Q{\ifmmode{{\mathbb Q}}\else{${\mathbb Q}$}\fi}
\def\C{\ifmmode{{\mathbb C}}\else{${\mathbb C}$}\fi}
\def\P{\ifmmode{{\mathbb P}}\else{${\mathbb P}$}\fi}

\def\H{\ifmmode{{\mathrm H}}\else{${\mathrm H}$}\fi}

\def\B{\ifmmode{{\cal B}}\else{${\cal B}$}\fi}
\def\E{\ifmmode{{\cal E}}\else{${\cal E}$}\fi}
\def\F{\ifmmode{{\cal F}}\else{${\cal F}$}\fi}
\def\K{\ifmmode{{\cal K}}\else{${\cal K}$}\fi}
\def\L{\ifmmode{{\cal L}}\else{${\cal L}$}\fi}
\def\M{\ifmmode{{\cal M}}\else{${\cal M}$}\fi}
\def\N{\ifmmode{{\cal N}}\else{${\cal N}$}\fi}
\def\O{\ifmmode{{\cal O}}\else{${\cal O}$}\fi}
\def\U{\ifmmode{{\cal U}}\else{${\cal U}$}\fi}
\def\X{\ifmmode{{\cal X}}\else{${\cal X}$}\fi}

\def\Br{\ifmmode{{\mathrm{Br}}}\else{${\mathrm{Br}}$}\fi}
\def\OG{\ifmmode{\widetilde{\cal M}_4}\else{$\widetilde{\cal M}_4$}\fi}
\def\D{\ifmmode{{\cal{D}}^b}\else{${{\cal{D}}^b}$}\fi}
\def\Shah{\ifmmode{\amalg\hspace*{-3.5pt}\amalg}\else{$\amalg\hspace*{-3.5pt}\amalg$}\fi}

\begin{document}

\title{Fourier-Mukai transforms, mirror symmetry, and generalized K3 surfaces}
\author{Justin Sawon}
\date{July, 2012}
\maketitle

\begin{abstract}

We study generalized complex structures on K3 surfaces, in the sense of Hitchin. For each real parameter $t\in (1,\infty)$ we exhibit two families of generalized K3 surfaces, $(M,\mathcal{I}_{\zeta})$ and $(M,\mathcal{J}_{\zeta})$, parametrized by $\zeta\in\P^1$, which are Mukai dual for $\zeta=0$ and $\infty$, and mirror partners for $\zeta\neq 0$ and $\infty$. Moreover, the Fourier-Mukai equivalence $D^b(M,\mathcal{I}_0)\rightarrow D^b(M,\mathcal{J}_0)$ induces an isomorphism $\phi_T$ between the spaces of first order deformations of $(M,\mathcal{I}_0)$ and $(M,\mathcal{J}_0)$ as generalized complex manifolds, and the deformations $(M,\mathcal{I}_{\zeta})$ and $(M,\mathcal{J}_{\zeta})$ agree under $\phi_T$, up to a $B$-field correction which vanishes in the limit $t\rightarrow\infty$.
\end{abstract}

\section{Introduction}

The results in this article were motivated by a desire to connect, in a concrete manner, Fourier-Mukai transforms and mirror symmetry. Now a Fourier-Mukai transform is an equivalence between the derived categories of two complex manifolds, whereas mirror symmetry relates complex and symplectic manifolds. Generalized complex geometry, introduced by Hitchin~\cite{hitchin03}, provides the framework for unifying complex and symplectic geometry; we can treat both complex and symplectic manifolds as examples of generalized complex manifolds. Moreover, it is possible to deform a complex manifold to a symplectic manifold in the category of generalized complex manifolds. In this article we start with a pair of Mukai dual (complex) K3 surfaces $X$ and $Y$ and deform them as generalized complex surfaces so that one remains complex while the other becomes a symplectic K3 surface. The resulting pair will be mirror partners.

The deformations of $X$ and $Y$ must be carefully chosen, and here we follow Toda~\cite{toda09}. First order deformations of the category $\mathrm{Coh}(X)$ of coherent sheaves on a complex manifold $X$ are parametrized by the degree two Hochschild cohomology $HH^2(X)$. This can be identified, via the HKR isomorphism, with
$$HT^2(X):=\H^2(X,\O_X)\oplus \H^1(X,T)\oplus \H^0(X,\wedge^2T).$$
Toda gave an explicit construction of the first order deformation of $\mathrm{Coh}(X)$ corresponding to an element $u\in HT^2(X)$. Now a Fourier-Mukai transform
$$\Phi:D^b(X)\rightarrow D^b(Y)$$
induces an isomorphism
$$\phi_T:HT^2(X)\rightarrow HT^2(Y).$$
Thus for every deformation of $\mathrm{Coh}(X)$, given by $u\in HT^2(X)$, there is a corresponding deformation of $\mathrm{Coh}(Y)$, given by $v:=\phi_T(u)\in HT^2(Y)$. Toda showed that the Fourier-Mukai transform $\Phi$ extends to an equivalence between the derived categories of these first order deformations.

Gualtieri~\cite{gualtieri11} showed that $HT^2(X)$ is also the space parametrizing first order deformations of the complex manifold $X$ as a generalized complex manifold. Our goal was to find a deformation of $X$ as a complex manifold, i.e., with
$$u\in\H^1(X,T)\subset HT^2(X),$$
such that the corresponding deformation $v\in HT^2(Y)$ takes $Y$ to a symplectic K3 surface. In this article, we achieve our goal in the limit. To be specific, we carefully choose the Mukai dual K3 surfaces $X$ and $Y$ so that for each real parameter $t\in (1,\infty)$ there is a deformation $u_t\in\H^1(X,T)$ of $X$ and a deformation $v_t\in HT^2(Y)$ taking $Y$ to a symplectic manifold, such that $v_t$ agrees with $\phi_T(u_t)$ up to a $B$-field correction. Moreover, the $B$-field correction vanishes as $t\rightarrow\infty$. The K3 surfaces are elliptic, and as one might expect, the limit $t\rightarrow\infty$ is the large complex structure limit where the fibres collapse.

Although Toda's results, the deformations of categories and Fourier-Mukai transforms, only apply to first order deformations, we are able to `integrate' our first order deformations $u_t$ and $v_t$ to produce genuine families of generalized K3 surfaces parametrized by $\zeta\in\P^1$. The deformation $(M,\mathcal{I}_{\zeta})$ of $X$ is a complex K3 surface, while for $\zeta\neq 0$ or $\infty$, the deformation $(M,\mathcal{J}_{\zeta})$ of $Y$ is a ($B$-field transform of a) symplectic K3 surface. We show that $(M,\mathcal{I}_{\zeta})$ and $(M,\mathcal{J}_{\zeta})$ are mirror partners, in the differential geometric sense of Gross~\cite{gross99}. This statement is true for all values of the parameter $t\in (1,\infty)$, not just in the limit.

One expects that $(M,\mathcal{I}_{\zeta})$ and $(M,\mathcal{J}_{\zeta})$ should be further related by Kontsevich's Homological Mirror Symmetry~\cite{kontsevich95}. Indeed, the required equivalence of categories should arise from a deformation of the original Fourier-Mukai transform
$$\Phi:D^b(X)\rightarrow D^b(Y).$$
However, a priori Toda's results only produce first order deformations, and the theory must be further developed before we can deal with global deformations.

The author would like to thank Andrei C{\u a}ld{\u a}raru, Marco Gaultieri, Mark Gross, and Emanuele Macr{\`i} for helpful explanations of their work.

\section{Statement of results}

The purpose of this article is to prove the following.

\begin{thm}
\label{main}
Let $X$ be an elliptic K3 surface which admits a section but is otherwise generic, let $Y$ be the Jacobian (i.e., dual elliptic fibration) of $X$, and let
$$\Phi:D^b(X)\rightarrow D^b(Y)$$
be the Fourier-Mukai transform between the derived categories of $X$ and $Y$ induced by the relative Poincar{\'e} line bundle. For each real parameter $t\in (1,\infty)$ there exists two families of generalized K3 surfaces, $(M,\mathcal{I}_{\zeta})$ and $(M,\mathcal{J}_{\zeta})$, parametrized by $\zeta\in\P^1=\mathbb{C}\cup\{\infty\}$, with the following properties:
\begin{enumerate}
\item $\mathcal{I}_{\zeta}$ is a complex structure for all $\zeta$, with $X=(M,\mathcal{I}_0)$ and $\overline{X}=(M,\mathcal{I}_{\infty})$,
\item $\mathcal{J}_{\zeta}$ is a complex structure for $\zeta=0$ and $\infty$, with $Y=(M,\mathcal{J}_0)$ and $\overline{Y}=(M,\mathcal{J}_{\infty})$, a symplectic structure for $|\zeta|=1$, and a $B$-field transform of a symplectic structure for all other values of $\zeta$,
\item to first order in $\zeta$, the deformation $(M,\mathcal{I}_{\zeta})$ of $X$ corresponds to the deformation $(M,\mathcal{J}_{\zeta})$ of $Y$ under the isomorphism
    $$\phi_T:HT^2(X)\rightarrow HT^2(Y)$$
    introduced by Toda~\cite{toda09}, up to a $B$-field correction which vanishes in the limit $t\rightarrow\infty$,
\item $(M,\mathcal{I}_{\zeta})$ and $(M,\mathcal{J}_{\zeta})$ are mirror K3 surfaces for $\zeta\neq 0$ and $\infty$, in the sense of Gross~\cite{gross99}.
\end{enumerate}
\end{thm}

\noindent
We start in Section~3 by describing the Fourier-Mukai transform $\Phi$ and the induced isomorphisms of the Hochschild structures. The two families will be constructed in Sections~4.2 and~4.3, respectively, and (1) and (2) of the theorem will then be immediate. Parts (3) and (4) will be clarified and proved in Sections~5.2 and~6.1, respectively. We consider the limiting behaviour as the real parameter $t\rightarrow 1^+$ and $t\rightarrow\infty$ in Section~5.3.

\section{Fourier-Mukai transforms}

\subsection{Mukai dual elliptic K3 surfaces}

Let $X$ be an elliptic K3 surface. Furthermore, assume that $X\rightarrow\P^1$ has a global section and that it is generic in the sense that there are $24$ nodal rational curves as singular fibres. Let $C_X$ and $F_X$ denote the section and a generic fibre of $X\rightarrow\P^1$; then $C^2=-2$, $F^2=0$, and $C.F=1$ (we will drop the subscripts when the notation is unambiguous).

Let $Y$ be the Jacobian of $X$, i.e., the dual elliptic fibration. Because $X\rightarrow\P^1$ has only nodal rational curves as singular fibres, $Y\rightarrow\P^1$ will be a smooth K3 surface. In fact, because $X\rightarrow\P^1$ admits a global section, $Y$ will be isomorphic to $X$ as an elliptic fibration; nevertheless, we will continue to distinguish $Y$ from $X$.

Since $Y$ parametrizes rank-one torsion-free sheaves on the fibres of $X\rightarrow\P^1$, there is a relative Poincar{\'e} line bundle $\mathcal{P}$ on the fibre product $X\times_{\P^1}Y$ (more accurately, $\mathcal{P}$ is a sheaf, as it may not be locally free at the nodes of the singular fibres). As a universal sheaf $\mathcal{P}$ is defined only up to tensoring with a line bundle pulled back from $Y$; the standard normalization requires both $\mathcal{P}|_{C_X\times_{\P^1}Y}$ and $\mathcal{P}|_{X\times_{\P^1}C_Y}$ to be trivial. We will use a different normalization, defining $\mathcal{L}:=\mathcal{P}\otimes\pi_{\P^1}^*\O(1)$ where $\pi_{\P^1}$ denotes the projection from $X\times_{\P^1}Y$ to $\P^1$ (cf.\ the functor $\mathbf{T}$ of Bartocci et al.~\cite{bbhm98}). This will produce a more symmetric cohomological Fourier-Mukai transform in Lemma~\ref{phiOmega}.

We can use $\mathcal{L}$ to construct a functor
\begin{eqnarray*}
\Phi: D^b(X) & \rightarrow & D^b(Y) \\
\mathcal{E}^{\bullet} & \mapsto & R\pi_{Y*}(\mathcal{L}\stackrel{L}{\otimes} \pi_X^*\mathcal{E}^{\bullet})
\end{eqnarray*}
where $\pi_X$ and $\pi_Y$ are the projections from $X\times_{\P^1}Y$ to $X$ and $Y$ respectively. It is well known that $\Phi$ will be a Fourier-Mukai transform in this case, i.e., an equivalence of triangulated categories (see Mukai~\cite{mukai87}).

\subsection{Hochschild structures}

We calculate here, for later use, the induced (co)homological Fourier-Mukai transforms. References for the following material are C{\u a}ld{\u a}raru~\cite{caldararu03i, caldararu05} and the revised version of the first article with Willerton~\cite{cw10}. Associated to a complex manifold $X$ is a graded ring $HH^{\bullet}(X)$ known as Hochschild cohomology and a graded $HH^{\bullet}(X)$-module, $HH_{\bullet}(X)$, known as Hochschild homology. The latter admits a non-degenerate pairing, the Mukai pairing
$$HH_{\bullet}(X)\times HH_{\bullet}(X)\rightarrow\C.$$
Although Hochschild (co)homology has nice functorial properties, for computations it is easier to work with the harmonic structure of $X$, which consists of the graded ring
$$HT^{\bullet}(X)=\bigoplus_{p+q=\bullet}\H^p(X,\Lambda^qT)$$
of cohomology of polyvector fields, and the graded $HT^{\bullet}(X)$-module
$$H\Omega_{\bullet}(X)=\bigoplus_{q-p=\bullet}\H^p(X,\Omega^q).$$
Note the unconventional choice of grading of the latter. There is also a Mukai pairing on $H\Omega_{\bullet}(X)$. The relation between the Hochschild structure and the harmonic structure is given by the Hochschild-Kostant-Rosenberg (HKR) isomorphisms:
$$I^{HKR}:HH^{\bullet}(X)\stackrel{\cong}{\longrightarrow}HT^{\bullet}(X)$$
$$I_{HKR}:HH_{\bullet}(X)\stackrel{\cong}{\longrightarrow}H\Omega_{\bullet}(X)$$
Note that these are isomorphisms of graded vector spaces. One difficulty is that they do not preserve all of the algebraic structure (i.e., the ring and module structures), but this can be corrected by twisting with the square root of the Todd polynomial. Thus we define compositions:
$$I^K:HH^{\bullet}(X)\stackrel{I^{HKR}}{\longrightarrow}HT^{\bullet}(X)\stackrel{\wedge Td_X^{-1/2}}{\longrightarrow}HT^{\bullet}(X)$$
$$I_K:HH_{\bullet}(X)\stackrel{I_{HKR}}{\longrightarrow}H\Omega_{\bullet}(X)\stackrel{\wedge Td_X^{1/2}}{\longrightarrow}H\Omega_{\bullet}(X)$$
We then have the following result, first conjectured by C{\u a}ld{\u a}raru~\cite{caldararu05}.

\begin{thm}[Calaque, Rossi, and Van den Bergh~\cite{crv09}]
\label{caldararu}
The maps $I^K$ and $I_K$ intertwine the module structures, i.e., the following diagram commutes, where the vertical arrows denote the respective actions:
$$\xymatrix@R=20pt@C=20pt{
	HH^{\bullet}(X)\times HH_{\bullet}(X)\ar[rr]^{I^K\times I_K}\ar[dd] & & HT^{\bullet}(X)\times H\Omega_{\bullet}(X) \ar[dd] \\
	 & & \\
	HH_{\bullet}(X)\ar[rr]^{I_K} & & H\Omega_{\bullet}(X) \\
	}$$
\end{thm}

\begin{rmk}
C{\u a}ld{\u a}raru's conjecture (now a theorem) is stronger than the above statement, as it also asserts that $I^K$ preserves the ring structures and $I_K$ preserves the Mukai pairings, but we will only need the above result. Moreover, we only need the result for K3 surfaces; there are earlier proofs in this case (in fact, for all Ricci-flat manifolds) by Huybrechts and Nieper-Wi{\ss}kirchen~\cite{hn11} and also by Macr{\`i}, Nieper-Wi{\ss}kirchen, and Stellari~\cite{mns08}.
\end{rmk}

Let $\Phi:D^b(X)\rightarrow D^b(Y)$ be a Fourier-Mukai transform with kernel $\U\in D^b(X\times Y)$. Since Hochschild (co)homology is functorial, $\Phi$ induces isomorphisms:
$$\phi:HH^{\bullet}(X)\stackrel{\cong}{\longrightarrow} HH^{\bullet}(Y)$$
$$\phi:HH_{\bullet}(X)\stackrel{\cong}{\longrightarrow} HH_{\bullet}(Y)$$
Again, we would like to transfer these to isomorphisms of the harmonic structure. To begin, we {\em define} an isomorphism $\phi_{HT}:HT^{\bullet}(X)\rightarrow HT^{\bullet}(Y)$ by $\phi_{HT}:=I^K\circ \phi\circ (I^K)^{-1}$. In other words, we define $\phi_{HT}$ so that the following diagram commutes:
$$\xymatrix@R=20pt@C=20pt{
	HH^{\bullet}(X)\ar[rr]^{\phi}\ar[dd]_{I^K} & & HH^{\bullet}(Y)\ar[dd]^{I^K} \\
	 & & \\
	HT^{\bullet}(X)\ar[rr]^{\phi_{HT}} & & HT^{\bullet}(Y) \\
	}$$
When it comes to $H\Omega_{\bullet}$, there is already a natural cohomological Fourier-Mukai transform
\begin{eqnarray*}
\phi_{H\Omega}:H\Omega_{\bullet}(X) & \rightarrow & H\Omega_{\bullet}(Y) \\
\alpha & \mapsto & \pi_{Y*}(v(\U)\wedge \pi_Y^*\alpha),
\end{eqnarray*}
where
$$v(\U)=ch(\U)\wedge Td_{X\times Y}^{1/2}\in H\Omega_0(X\times Y)$$
is the Mukai vector of $\U$. Fortunately, this is also compatible with the isomorphism of Hochschild homology.

\begin{lem}
\label{commute1}
The following diagram commutes:
$$\xymatrix@R=20pt@C=20pt{
	HH_{\bullet}(X)\ar[rr]^{\phi}\ar[dd]_{I_K} & & HH_{\bullet}(Y)\ar[dd]^{I_K} \\
	 & & \\
	H\Omega_{\bullet}(X)\ar[rr]^{\phi_{H\Omega}} & & H\Omega_{\bullet}(Y) \\
	}$$
\end{lem}

\begin{prf}
This was proved by Macr{\'i} and Stellari~\cite[Theorem~1.2]{ms09}, building on results of Ramadoss~\cite{ramadoss08}.
\end{prf}

Assume now that $X$ and $Y$ are K3 surfaces, and denote by
$$\sigma_X\in H\Omega_2(X)=\H^0(X,\Omega^2)$$
a holomorphic two-form on $X$, and by $\sigma_Y:=\phi_{H\Omega}(\sigma_X)$ the corresponding two-form on $Y$. There are also corresponding elements of Hochschild homology $HH_2(X)$ and $HH_2(Y)$, given by $I_K^{-1}(\sigma_X)$ and $I_K^{-1}(\sigma_Y)$, respectively. By Lemma~\ref{commute1} we have
$$\phi(I_K^{-1}(\sigma_X))=I_K^{-1}(\phi_{H\Omega}(\sigma_X))=I_K^{-1}(\sigma_Y).$$
The action of $HH^2(X)$ on $I_K^{-1}(\sigma_X)\in HH_2(X)$ induces an isomorphism
$$\neg I_K^{-1}(\sigma_X):HH^2(X)\longrightarrow HH_0(X),$$
and likewise for $Y$. By functoriality, the following diagram commutes:
$$\xymatrix@R=20pt@C=20pt{
	HH^2(X)\ar[rr]^{\phi}\ar[dd]_{\neg I_K^{-1}(\sigma_X)} & & HH^2(Y)\ar[dd]^{\neg I_K^{-1}(\sigma_Y)} \\
	 & & \\
	HH_0(X)\ar[rr]^{\phi} & & HH_0(Y) \\
	}$$
Similarly, the action of $HT^2(X)$ on $\sigma_X\in H\Omega_2(X)$ induces an isomorphism
$$\neg \sigma_X:HT^2(X)\longrightarrow H\Omega_0(X),$$
and likewise for $Y$. Putting everything together gives the next result.

\begin{prp}
\label{phiprime}
The following diagram commutes:
$$\xymatrix@R=20pt@C=20pt{
	HT^2(X)\ar[rr]^{\phi_{HT}}\ar[dd]_{\neg\sigma_X} & & HT^2(Y)\ar[dd]^{\neg\sigma_Y} \\
	 & & \\
	H\Omega_0(X)\ar[rr]^{\phi_{H\Omega}} & & H\Omega_0(Y) \\
	}$$
\end{prp}

\begin{prf}
Consider the following diagram:
$$\xymatrix@R=25pt@C=25pt{
	& HT^2(X) \ar[rr]^(0.6){\phi_{HT}}\ar'[d][dd]_{\neg\sigma_X}
	& & HT^2(Y) \ar[dd]^(0.7){\neg\sigma_Y}
	\\
	HH^2(X) \ar[ur]^{I^K}\ar[rr]^(0.6){\phi}\ar[dd]_(0.65){\neg I_K^{-1}(\sigma_X)}
	& & HH^2(Y) \ar[ur]^{I^K}\ar[dd]^(0.65){\neg I_K^{-1}(\sigma_Y)}
	\\
	& H\Omega_0(X) \ar'[r][rr]^{\hspace*{-12mm}\phi_{H\Omega}}
	& & H\Omega_0(Y)
	\\
	HH_0(X) \ar[rr]^(0.6){\phi}\ar[ur]_{I_K}
	& & HH_0(Y) \ar[ur]_{I_K}
	}$$
The left and right faces commute by Theorem~\ref{caldararu}. The top face commutes by the definition of $\phi_{HT}$. The bottom face commutes by Lemma~\ref{commute1}. The front face commutes by functoriality of Hochschild cohomology and homology, as mentioned above. It follows that the whole diagram commutes; in particular, the back face commutes, which proves the proposition.
\end{prf}

\begin{rmk}
The point is that $\phi_{H\Omega}$ can be calculated directly, and then the above diagram can be used to determine $\phi_{HT}$. We actually wish to determine the isomorphism
$$\phi_T:HT^2(X)\longrightarrow HT^2(Y)$$
{\em defined} by Toda~\cite{toda09} as $\phi_T:=I^{HKR}\circ\phi\circ (I^{HKR})^{-1}$. This differs from $\phi_{HT}$ insofar as the definition uses $I^{HKR}$ instead of $I^K$; therefore
$$\phi_T=(\wedge Td_Y^{1/2})\circ \phi_{HT}\circ (\wedge Td_X^{-1/2}).$$
\end{rmk}

We return now to our pair of dual elliptic K3 surfaces, $X\rightarrow\P^1$ and $Y\rightarrow\P^1$. Denote by $1_X$ and $\eta_X$ the generators of $\H^0(X,\Z)$ and $\H^4(X,\Z)$ respectively, and by $[C_X]$ and $[F_X]$ the classes of the section and a fibre in $\H^2(X,\Z)\cap\H^{1,1}(X)$. These can be regarded as elements in
$$H\Omega_0(X)=\H^0(X,\Omega^0)\oplus\H^1(X,\Omega^1)\oplus\H^2(X,\Omega^2).$$
Denote by $1_Y$, $\eta_Y$, $[C_Y]$, and $[F_Y]$ the corresponding elements of $H\Omega_0(Y)$.

\begin{lem}
\label{phiOmega}
The cohomological Fourier-Mukai transform
$$\phi_{H\Omega} : H\Omega_0(X)\rightarrow H\Omega_0(Y)$$
is given by
\begin{eqnarray*}
1_X & \mapsto & -[C_Y]-[F_Y] \\
\eta_X & \mapsto & [F_Y] \\
{[C_X]} & \mapsto & 1_Y+\eta_Y \\
{[F_X]} & \mapsto & -\eta_Y.
\end{eqnarray*}
\end{lem}

\begin{prf}
This is essentially Lemma~4.1 of Bartocci et al.~\cite{bbhm98} or Lemma~7.11 of Huybrechts~\cite{huybrechts04}, but they use slightly different normalizations so we give a complete proof here. We use Theorem~1.1 of Mukai~\cite{mukai87}, which says that if $\hat{\Phi}: D^b(Y)\rightarrow D^b(X)$ is the Fourier-Mukai transform with the same kernel $\mathcal{L}$ but in the opposite direction, then there are natural isomorphisms
$$\hat{\Phi}\circ\Phi\cong(-\mathrm{Id}_X)^*[-1],$$
$$\Phi\circ\hat{\Phi}\cong(-\mathrm{Id}_Y)^*[-1],$$
where $-\mathrm{Id}_X$ is the map given by multiplication by $-1$ on the fibres (and similarly for $Y$) and $[-1]$ is the shift functor. Note that Mukai's Theorem~1.1 is stated for the Fourier-Mukai transforms with kernel $\mathcal{P}$, and so it includes $\otimes\omega^{-1}_{X/\P^1}$ on the right hand side of the first isomorphism, where $\omega_{X/\P^1}\cong\pi_{\P^1}^*\O(2)$ is the relative dualizing sheaf of $X\rightarrow\P^1$ (and similarly for the second isomorphism). These terms are eliminated by our choice of normalization.

Observe that $\Phi$ takes the skyscraper sheaf $\O_{F_X\cap C_X}$ supported at the point $F_X\cap C_X$ to the structure sheaf $\O_{F_Y}$ of the corresponding fibre of $Y\rightarrow\P^1$. The above result of Mukai (or a direct observation) shows that $\Phi$ also takes the structure sheaf $\O_{F_X}$ of a fibre of $X\rightarrow\P^1$ to $\O_{F_Y\cap C_Y}[-1]$, the shifted skyscraper sheaf supported at $F_Y\cap C_Y$. Taking Mukai vectors yields $\eta_X\mapsto [F_Y]$ and $[F_X]\mapsto -\eta_Y$.

Next observe that $\Phi$ takes $\O_X$ to $\O_{C_Y}(-1)[-1]$ (cf.~\cite[Theorem~2.8(4)]{bbhm98}, noting the different choice of normalization). Equivalently, by Mukai's result, $\Phi$ takes $\O_{C_X}$ to $\O_Y(1)$. Taking Mukai vectors yields $1_X+\eta_X\mapsto -[C_Y]$ and $[C_X]+\eta_X\mapsto 1_Y+[F_Y]+\eta_Y$. The lemma follows.
\end{prf}

\begin{rmk}
Note that $\phi_{H\Omega}$ must be an isometry, i.e., it must preserve the Mukai pairing, which is easily verified in this case.
\end{rmk}

Next we will calculate $\phi_{HT}$. Write
$$HT^2(X)=\H^0(X,\Lambda^2T)\oplus\H^1(X,T)\oplus\H^2(X,\O),$$
with $\H^0(X,\Lambda^2T)$ and $\H^2(X,\O)$ generated by $\sigma_X^{-1}$ and $\bar{\sigma}_X$, respectively.

\begin{lem}
\label{isomorphism}
The map
$$\neg\sigma_X:HT^2(X)\longrightarrow H\Omega_0(X)$$
is given by
\begin{eqnarray*}
\sigma^{-1}_X & \mapsto & 4.1_X \\
\bar{\sigma}_X & \mapsto & \sigma_X\bar{\sigma}_X \\
\sigma^{-1}_X[C_X] & \mapsto & [C_X] \\
\sigma^{-1}_X[F_X] & \mapsto & [F_X].
\end{eqnarray*}
\end{lem}

\begin{prf}
The first line follows from the local calculation
$$\left(\sigma^{ij}\frac{\partial\phantom{z}}{\partial z^i}\wedge\frac{\partial\phantom{z}}{\partial z^j}\right)\neg \left(\sigma_{kl}dz^k\wedge dz^l\right)=\sigma^{ij}\sigma_{kl}\left(\delta^k_i\delta^l_j-\delta^k_j\delta^l_i\right)=2\sigma^{ij}\sigma_{ij}=4.$$
The second is vacuous.

Note that $\sigma_X$ induces a bundle isomorphism $T\cong\Omega^1$, given by $w\mapsto\sigma(w,-)$, which in turn induces an isomorphism $\H^1(X,T)\cong\H^1(X,\Omega^1)$. Then we denote by $
\sigma^{-1}_X[C_X]$ and $\sigma^{-1}_X[F_X]$ the elements of $\H^1(X,T)$ corresponding to $[C_X]$ and $[F_X]\in \H^1(X,\Omega^1)$, respectively. The third and fourth lines are then automatic by definition.
\end{prf}

\begin{rmk}
Since $\H^0(X,\Omega^2)\cong\C$, the holomorphic two-form $\sigma_X$ is determined up to scale; we will (partially) normalize by requiring
$$\sigma_X\bar{\sigma}_X=4\eta_X.$$
In the presence of a hyperk{\"a}hler structure, $\sigma_X=\omega_J+i\omega_K$ and
$$\int_X\sigma_X\bar{\sigma}_X=\int_X\omega_J^2+\omega_K^2=4\int_X\frac{\omega_I^2}{2}=4\mathrm{vol}(X),$$
so our choice of normalization corresponds to $\mathrm{vol}(X)=1$. Note that
$$\int_Y\sigma_Y\bar{\sigma}_Y=\int_X\sigma_X\bar{\sigma}_X$$
because $\sigma_Y=\phi_{H\Omega}(\sigma_X)$ and $\phi_{H\Omega}$ preserves the Mukai pairing, therefore we also have
$$\sigma_Y\bar{\sigma}_Y=4\eta_Y.$$
\end{rmk}

\begin{lem}
\label{phi'_T}
The isomorphism
$$\phi_{HT} : HT^2(X)\longrightarrow HT^2(Y)$$
is given by
\begin{eqnarray*}
\frac{1}{4}\sigma_X^{-1} & \mapsto & -\sigma_Y^{-1}[C_Y]-\sigma_Y^{-1}[F_Y] \\
\frac{1}{4}\bar{\sigma}_X & \mapsto & \sigma_Y^{-1}[F_Y] \\
\sigma_X^{-1}[C_X] & \mapsto & \frac{1}{4}\sigma_Y^{-1}+\frac{1}{4}\bar{\sigma}_Y \\
\sigma_X^{-1}[F_X] & \mapsto & -\frac{1}{4}\bar{\sigma}_Y. \\
\end{eqnarray*}
\end{lem}

\begin{prf}
This follows from Proposition~\ref{phiprime} and Lemmas~\ref{phiOmega} and~\ref{isomorphism}.
\end{prf}

\begin{lem}
\label{phi_T}
The isomorphism
$$\phi_T : HT^2(X)\longrightarrow HT^2(Y)$$
is given by
\begin{eqnarray*}
\frac{1}{4}\sigma_X^{-1} & \mapsto & -\sigma_Y^{-1}[C_Y]-2\sigma_Y^{-1}[F_Y] \\
\frac{1}{4}\bar{\sigma}_X & \mapsto & \sigma_Y^{-1}[F_Y] \\
\sigma_X^{-1}[C_X] & \mapsto & \frac{1}{4}\sigma_Y^{-1}+\frac{2}{4}\bar{\sigma}_Y \\
\sigma_X^{-1}[F_X] & \mapsto & -\frac{1}{4}\bar{\sigma}_Y. \\
\end{eqnarray*}
\end{lem}

\begin{prf}
Recall that
$$\phi_T=(\wedge Td_Y^{1/2})\circ \phi_{HT}\circ (\wedge Td_X^{-1/2}).$$
So we simply compose $\phi_{HT}$ with $Td_X^{-1/2}=1_X-\eta_X$ and $Td_Y^{1/2}=1_Y+\eta_Y$ to get the result.
\end{prf}

\section{The construction and basic properties}

\subsection{Generalized K3 surfaces}

Generalized complex geometry was introduced by Hitchin~\cite{hitchin03}, and subsequently developed by Gualtieri~\cite{gualtieri11}. Let $M$ be a smooth manifold with tangent bundle $T$. A generalized complex structure on $M$ is an endomorphism $\mathcal{J}$ of $T\oplus T^*$ such that $\mathcal{J}^2=-\mathrm{Id}$, $\mathcal{J}$ is orthogonal with respect to the natural inner product on $T\oplus T^*$, and $\mathcal{J}$ satisfies a certain integrability condition. Complex manifolds and symplectic manifolds both give examples of generalized complex manifolds: a complex structure $I$ and a symplectic structure $\omega$ induce generalized complex structures
$$\mathcal{J}_I:=\left(\begin{array}{cc} -I & 0 \\ 0 & I^* \end{array}\right)\qquad\mbox{and}\qquad\mathcal{J}_{\omega}:=\left(\begin{array}{cc} 0 & -\omega^{-1} \\ \omega & 0 \end{array}\right)$$
respectively. A generalized complex structure is determined by a pure spinor line, which is given locally by a section of $\wedge^{\bullet}T^*\otimes\C$. The pure spinors for $\mathcal{J}_I$ and $\mathcal{J}_{\omega}$ above are
$$dz_1\wedge\ldots\wedge dz_n\in\Omega^{n,0}\qquad\mbox{and}\qquad\mathrm{exp}(i\omega)\in\wedge^{\mathrm{even}}T^*\otimes\C,$$
respectively.

\begin{dfn}
Let $M$ be the underlying smooth $4$-manifold of a K3 surface. In this article, by a generalized K3 surface we shall mean $M$ equipped with a generalized complex structure.
\end{dfn}

\subsection{The first family}

Let $X$ be a (complex) K3 surface, and let $\alpha$ be a K{\"a}hler class on $X$. By Yau's Theorem there exists a hyperk{\"a}hler metric $g$ on $X$ whose K{\"a}hler form $\omega_I$ represents $\alpha$. Let $\sigma_X$ be a holomorphic two-form on $X$, normalized so that
$$\int_X\sigma_X\bar{\sigma}_X=2\int_X\omega_I^2=4\mathrm{vol}(X,g),$$
and write $\omega_J$ and $\omega_K$ for the real and imaginary parts of $\sigma$.

We recall twistor families of K3 surfaces, following Section 3(F) of Hitchin et al.~\cite{hklr87}. If $M$ is the underlying $4$-manifold of $X$, then there is a family of complex structures $aI+bJ+cK$, parametrized by $(a,b,c)\in S^2$, which are compatible with $g$, i.e., which make $(M,g)$ into a K{\"a}hler manifold. Change to a single complex parameter $\zeta\in\P^1=\mathbb{C}\cup\{\infty\}$ and write the complex structure as $I_{\zeta}$. Then $I_{\zeta}$ is determined by the holomorphic two-form on $(M,I_{\zeta})$, which up to scale is given by
$$\sigma_{\zeta}:=\sigma_X+2\zeta\omega_I-\zeta^2\bar{\sigma}_X.$$
Note that $X=(M,I_0)$, while $I_{\infty}=-I$ gives the conjugate complex surface $\overline{X}$. The generic K3 surface $(M,I_{\zeta})$ in this twistor family is non-algebraic.

\begin{prf}{\bf of Theorem~\ref{main}(1)}
Consider our elliptic K3 surface $X$, with section $C$ and fibre $F$. If we define
$$\alpha:=\frac{1}{t}[C]+\left(\frac{t^2+1}{t}\right)[F]\in\H^{1,1}(X)\cap\H^2(X,\mathbb{R})$$
with $t\in\mathbb{R}$, then
$$\alpha.C=\frac{1}{t}C^2+\left(\frac{t^2+1}{t}\right)F.C=\frac{t^2-1}{t},\qquad \alpha.F=\frac{1}{t}C.F+\left(\frac{t^2+1}{t}\right)F^2=\frac{1}{t},$$
$$\mbox{and}\qquad \alpha^2=\frac{1}{t^2}C^2+2\left(\frac{t^2+1}{t^2}\right)C.F+\left(\frac{t^2+1}{t}\right)^2F^2=2.$$
So $\alpha$ is a K{\"a}hler class on $X$ provided $t>1$ (note that the N{\'e}ron-Severi group $NS(X)$ is generated by $C$ and $F$ because $X$ is generic). As above, there exists a hyperk{\"a}hler metric $g$ on $X$ whose K{\"a}hler form $\omega_I$ represents $\alpha$, and a corresponding twistor family $(M,I_{\zeta})$. We regard this simply as a family of complex K3 surfaces; in other words, we ignore the metric and K{\"a}hler structures. Note that
$$\int_X\sigma_X\bar{\sigma}_X=2\int_X\omega_I^2=2\int_X\alpha^2=4,$$
so our normalization of $\sigma_X$ agrees with the one used in Section~3.

Finally, we define $\mathcal{I}_{\zeta}$ to be the generalized complex structure $\mathcal{I}_{I_{\zeta}}$ coming from $I_{\zeta}$. The first family of generalized K3 surfaces, $(M,\mathcal{I}_{\zeta})$, clearly satisfies the requirements of Theorem~\ref{main}(1).
\end{prf}

\subsection{The second family}

Let $Y$ be a (complex) K3 surface, with complex structure $I$. Choose a holomorphic two-form $\sigma$ on $Y$ (unique up to scale), and let $\omega_J$ and $\omega_K$ be the real and imaginary parts of $\sigma$. Then $\omega_J$ and $\omega_K$ are symplectic forms on the underlying $4$-manifold $M$. Gualtieri~\cite{gualtieri11} showed that we can interpolate between the generalized complex structure $\mathcal{J}_I$ of complex type and the generalized complex structure $\mathcal{J}_{\omega_J}$ of symplectic type.

\begin{lem}
Define
\begin{eqnarray*}
\mathcal{J}_{\theta} & := & (\cos\theta)\mathcal{J}_I+(\sin\theta)\mathcal{J}_{\omega_J} \\
 & = & \cos\theta\left(\begin{array}{cc} -I & 0 \\ 0 & I^* \\ \end{array}\right)+\sin\theta\left(\begin{array}{cc} 0 & -\omega^{-1}_J \\ \omega_J & 0 \\ \end{array}\right)
\end{eqnarray*}
for $\theta\in\mathbb{R}$. If $\theta$ is not an integer multiple of $\pi$, then $\mathcal{J}_{\theta}$ is a $B$-field transform of $\mathcal{J}_{(\csc\theta)\omega_J}$. Consequently, $\mathcal{J}_{\theta}$ is a generalized complex structure for all $\theta$.
\end{lem}

\begin{prf}
This is essentially Proposition 3.31 of~\cite{gualtieri11}. A direct calculation shows that $\mathcal{J}_{\theta}^2=-\mathrm{Id}$ and that
$$\mathcal{J}_{\theta}=e^{-B}\mathcal{J}_{(\csc\theta)\omega_J}e^B=\left(\begin{array}{cc} 1 & 0 \\ -B & 1 \\ \end{array}\right)
\left(\begin{array}{cc} 0 & -((\csc\theta)\omega_J)^{-1} \\ (\csc\theta)\omega_J & 0 \\ \end{array}\right)\left(\begin{array}{cc} 1 & 0 \\ B & 1 \\ \end{array}\right)$$
where $B=-(\cot\theta)\omega_K$. Since $B$-field transforms preserve integrability, $\mathcal{J}_{\theta}$ is integrable if $\theta$ is not an integer multiple of $\pi$. On the other hand, $\mathcal{J}_{2k\pi}=\mathcal{J}_I$ and $\mathcal{J}_{(2k+1)\pi}=-\mathcal{J}_I=\mathcal{J}_{-I}$ are also integrable.
\end{prf}

\noindent
We can replace $\omega_J$ by a combination $(\cos\phi)\omega_J+(\sin\phi)\omega_K$, which corresponds to multiplying the holomorphic two-form $\sigma$ by $e^{-i\phi}$.

\begin{dfn}
Define
$$\mathcal{J}_{\theta,\phi}:=(\cos\theta)\mathcal{J}_I+(\sin\theta\cos\phi)\mathcal{J}_{\omega_J}+(\sin\theta\sin\phi)\mathcal{J}_{\omega_K}.$$
By the above lemma, $\mathcal{J}_{\theta,\phi}$ is a generalized complex structure for all $\theta$ and $\phi\in\mathbb{R}$. Moreover, $\theta$ and $\phi$ may be thought of as spherical coordinates on the sphere $S^2$; changing to a single complex parameter $\zeta\in\P^1=\mathbb{C}\cup\{\infty\}$ gives
$$\mathcal{J}_{\zeta}=\left(\frac{1-|\zeta|^2}{1+|\zeta|^2}\right)\mathcal{J}_I+\left(\frac{-2\mathcal{I}m\zeta}{1+|\zeta|^2}\right)\mathcal{J}_{\omega_J}+\left(\frac{2\mathcal{R}e\zeta}{1+|\zeta|^2}\right)\mathcal{J}_{\omega_K}.$$
\end{dfn}

\begin{rmk}
Straight-up stereographic projection of $S^2$ onto the complex plane gives
$$(\cos\theta,\sin\theta\cos\phi,\sin\theta\sin\phi)=\left(\frac{1-|\zeta|^2}{1+|\zeta|^2},\frac{2\mathcal{R}e\zeta}{1+|\zeta|^2},\frac{2\mathcal{I}m\zeta}{1+|\zeta|^2}\right).$$
For our choice of complex parametrization we have also rotated by $\pi/2$ in the $\phi$ direction, or equivalently, multiplied the complex parameter by $i$. This is mainly done to improve the appearance of later formulae, but it can be justified by the fact that going from $\zeta=0$ to $1$ now represents a deformation from $\mathcal{J}_I$ to $\mathcal{J}_{\omega_K}$, which is the traditional hyperk{\"a}hler rotation, i.e., a holomorphic Lagrangian fibration with respect to $I$ becomes a special Lagrangian fibration with respect to $\omega_K$.
\end{rmk}

\begin{lem}
\label{spinor}
The pure spinor defining $\mathcal{J}_{\zeta}$ is given by
$$\sigma+2\zeta\left(1-\frac{1}{4}\sigma\bar{\sigma}\right)-\zeta^2\bar{\sigma}.$$
In particular, the family of generalized complex structures $\mathcal{J}_{\zeta}$ depends holomorphically on $\zeta\in\P^1$.
\end{lem}

\begin{prf}
For $\theta$ not an integer multiple of $\pi$, $\mathcal{J}_{\theta,\phi}$ is a $B$-field transform of the generalized complex structure coming from the symplectic structure $\omega=\csc\theta((\cos\phi)\omega_J+(\sin\phi)\omega_K)$, where $B=\cot\theta((\sin\phi)\omega_J-(\cos\phi)\omega_K)$. The relations between spherical coordinates $(\theta,\phi)$ and the complex coordinate $\zeta$ give
$$\cos\theta=\frac{1-|\zeta|^2}{1+|\zeta|^2}\qquad\mbox{and}\qquad\sin\theta=\frac{2|\zeta|}{1+|\zeta|^2}.$$
We can therefore write
$$B=\frac{1-|\zeta|^2}{2|\zeta|}\mathcal{I}m(-e^{-i\phi}\sigma)=\frac{1-|\zeta|^2}{2|\zeta|^2}\mathcal{R}e(\bar{\zeta}\sigma),$$
$$\omega=\frac{1+|\zeta|^2}{2|\zeta|}\mathcal{R}e(e^{-i\phi}\sigma)=\frac{1+|\zeta|^2}{2|\zeta|^2}\mathcal{I}m(\bar{\zeta}\sigma),$$
and
\begin{eqnarray*}
B+i\omega & = & \frac{1-|\zeta|^2}{2|\zeta|^2}\mathcal{R}e(\bar{\zeta}\sigma)+i\frac{1+|\zeta|^2}{2|\zeta|^2}\mathcal{I}m(\bar{\zeta}\sigma) \\
 & = & \frac{1}{2|\zeta|^2}(\mathcal{R}e(\bar{\zeta}\sigma)+i\mathcal{I}m(\bar{\zeta}\sigma))-\frac{1}{2}(\mathcal{R}e(\bar{\zeta}\sigma)-i\mathcal{I}m(\bar{\zeta}\sigma)) \\
 & = & \frac{\sigma}{2\zeta}-\frac{\zeta\bar{\sigma}}{2}. \\
 \end{eqnarray*}
The corresponding pure spinor is thus
\begin{eqnarray*}
e^Be^{i\omega} & = & \exp\left(\frac{\sigma}{2\zeta}-\frac{\zeta\bar{\sigma}}{2}\right) \\
 & = & 1+\left(\frac{\sigma}{2\zeta}-\frac{\zeta\bar{\sigma}}{2}\right)-\frac{1}{4}\sigma\bar{\sigma}. \\
\end{eqnarray*}
Since the pure spinor is only defined up to scale, we can multiply by $2\zeta$ to yield the desired result. Notice that the formula also gives the correct pure spinors, $\sigma$ and $\bar{\sigma}$ for $\mathcal{J}_I$ and $\mathcal{J}_{-I}$ respectively, when $\theta$ is a multiple of $\pi$, i.e., when $\zeta=0$ or $\infty$.
\end{prf}

These arguments apply to any K3 surface. Notice that we did not specify a K{\"a}her form or metric; the construction of this family depends only on the complex structure and choice of a holomorphic two-form on $Y$.

\begin{prf}{\bf of Theorem~\ref{main}(2)}
Consider the Jacobian $Y$ of the elliptic K3 surface $X$ introduced earlier. As we saw in Section~3, the Fourier-Mukai transform
$$\sigma_Y:=\phi_{H\Omega}(\sigma_X)$$
of the holomorphic two-form on $X$ gives a holomorphic two-form on $Y$, normalized so that
$$\int_Y\sigma_Y\bar{\sigma}_Y=4.$$
However, we will rescale and use $\sigma=t\sigma_Y$ instead. Thus we define $\mathcal{J}_{\zeta}$ to be the family of generalized complex structure associated to $Y$ and $t\sigma_Y$ as described above. The second family of generalized K3 surfaces, $(M,\mathcal{J}_{\zeta})$, then satisfies the required properties of Theorem~\ref{main}(2):
\begin{itemize}
\item $\mathcal{J}_0=\mathcal{J}_I$ is a complex structure, and $(M,\mathcal{J}_0)$ is the original K3 surface $Y=(M,I)$,
\item $\mathcal{J}_{\infty}=\mathcal{J}_{-I}$ is also a complex structure, and $(M,\mathcal{J}_{\infty})$ is the conjugate complex surface $\overline{Y}=(M,-I)$,
\item when $\zeta=e^{i(\phi-\pi/2)}$, i.e., when $|\zeta|=1$, $\mathcal{J}_{e^{i(\phi-\pi/2)}}=\mathcal{J}_{(\cos\phi)\omega_J+(\sin\phi)\omega_K}$ is a symplectic structure,
\item for all other values of $\zeta$, $\mathcal{J}_{\zeta}$ is a $B$-field transform of a symplectic structure.
\end{itemize}
\end{prf}

\begin{rmk}
Note that we now have
$$B+i\omega=\frac{t\sigma_Y}{2\zeta}-\frac{\zeta t\bar{\sigma}_Y}{2}$$
and the pure spinor is
$$t\sigma_Y+2\zeta\left(1-\frac{1}{4}t^2\sigma_Y\bar{\sigma}_Y\right)-\zeta^2t\bar{\sigma}_Y,$$
which we can again rescale to
$$\sigma_Y+2\zeta\left(\frac{1}{t}-\frac{1}{4}t\sigma_Y\bar{\sigma}_Y\right)-\zeta^2\bar{\sigma}_Y.$$
\end{rmk}

\section{Deformation theory}

\subsection{Deformations of generalized complex manifolds}

Gualtieri~\cite{gualtieri11} developed the deformation theory of generalized complex manifolds. Let $(M,\mathcal{J})$ be a generalized complex manifold and denote by $L$ the $+i$-eigenspace of $\mathcal{J}$ in $(T\oplus T^*)\otimes\C$. Then $L$ is a Lie algebroid and there exists an elliptic differential complex
$$(C^{\infty}(\Lambda^{\bullet}L^*),d_L).$$
The first order deformations of the generalized complex structure are parametrized by the degree two cohomology $\H^2(M,L)$; obstructions to deforming lie in $\H^3(M,L)$. If we regard a complex manifold $X=(M,I)$ as a generalized complex manifold, with $\mathcal{J}=\mathcal{J}_I$, then $L$ is $T^{0,1}\oplus\Omega^{1,0}$, $L^*$ is $\Omega^{0,1}\oplus T^{1,0}$, $d_L$ is $\bar{\partial}$, and $\H^2(M,L)$ is
$$HT^2(X):=\H^0(X,\Lambda^2T)\oplus\H^1(X,T)\oplus\H^2(X,\mathcal{O}).$$

In Section~4 we introduced two families of generalized K3 surfaces, $(M,\mathcal{I}_{\zeta})$ and $(M,\mathcal{J}_{\zeta})$, which are complex surfaces $X=(M,\mathcal{I}_0)$ and $Y=(M,\mathcal{J}_0)$ when $\zeta=0$. In this section we will determine the directions of these deformations at $\zeta=0$, i.e., we will calculate their classes in $HT^2(X)$ and $HT^2(Y)$. The obstructions must vanish since these deformations come from actual families; in any case, $HT^3(X)$ vanishes for K3 surfaces. More generally, Goto~\cite{goto05} has proved the unobstructedness of deformations of Calabi-Yau manifolds as generalized complex manifolds.

\begin{lem}
\label{direction1}
The family $(M,\mathcal{I}_{\zeta})$ is a deformation of $X=(M,\mathcal{I}_0)$ in the direction
$$-\frac{2}{t}\sigma_X^{-1}[C_X]-\frac{2(t^2+1)}{t}\sigma_X^{-1}[F_X]\in\H^1(X,T)\subset HT^2(X).$$
\end{lem}

\begin{prf}
Recall that $(M,\mathcal{I}_{\zeta})$ is a family of complex K3 surfaces. We denote the complex structures by $I_{\zeta}$ and write $X_{\zeta}:=(M,I_{\zeta})$. Up to scale, the holomorphic two-form on $X_{\zeta}$ is given by
$$\sigma_{\zeta}=\sigma_X+2\zeta\omega_I-\zeta^2\bar{\sigma}_X,$$
where the K{\"a}hler form $\omega_I$ represents the class $\alpha=\frac{1}{t}[C_X]+\left(\frac{t^2+1}{t}\right)[F_X]$.

Let $Z\in T\otimes\mathbb{C}$ and write $Z=Z^{1,0}+Z^{0,1}$ for the decomposition with respect to $I_0$. Then $Z$ is of type $(0,1)$ with respect to $I_{\zeta}$ if and only if $\sigma_{\zeta}(Z,-)=0$, because $\sigma_{\zeta}$ is a $(2,0)$-form on $X_{\zeta}$. In particular, the $(1,0)$-component of the one-form $\sigma_{\zeta}(Z,-)$ must vanish, and this is given by
$$\sigma_X(Z^{1,0},-)+2\zeta\omega_I(Z^{0,1},-).$$
Therefore
$$Z^{1,0}=-2\zeta\sigma_X^{-1}(\omega_I(Z^{0,1},-)),$$
and we see that to first order in $\zeta$ the vectors of type $(0,1)$ on $X_{\zeta}$ lie on the graph of the map $T^{0,1}\rightarrow T^{1,0}$ given by the element
$$-2\zeta\sigma_X^{-1}\omega_I\in\Omega^{0,1}\otimes T^{1,0}.$$
This means that we are deforming in the direction
$$-2\sigma_X^{-1}[\omega_I]=-\frac{2}{t}\sigma_X^{-1}[C_X]-\frac{2(t^2+1)}{t}\sigma_X^{-1}[F_X]\in\H^{0,1}_{\bar{\partial}}(X,T).$$
\end{prf}

\begin{rmk}
It is usually stated that a twistor family is a deformation in the direction $\sigma_X^{-1}[\omega_I]$, but this would correspond to the family of holomorphic two-forms
$$\sigma_X-\zeta\omega_I-\frac{\zeta^2}{4}\bar{\sigma}_X.$$
Of course, this just comes from a different parametrization of $\P^1$.
\end{rmk}

\begin{lem}
\label{direction2}
The family $(M,\mathcal{J}_{\zeta})$ is a deformation of $Y=(M,\mathcal{J}_0)$ in the direction
$$\frac{1}{2}\left(-\frac{1}{t}\sigma_Y^{-1}+t\bar{\sigma}_Y\right)\in\H^0(Y,\Lambda^2T)\oplus\H^2(Y,\mathcal{O})\subset HT^2(Y).$$
\end{lem}

\begin{prf}
Recall that the family of generalized complex structures is given by
$$\mathcal{J}_{\zeta}:=\left(\frac{1-|\zeta|^2}{1+|\zeta|^2}\right)\mathcal{J}_I+\left(\frac{-2\mathcal{I}m\zeta}{1+|\zeta|^2}\right)\mathcal{J}_{\omega_J}+\left(\frac{2\mathcal{R}e\zeta}{1+|\zeta|^2}\right)\mathcal{J}_{\omega_K}.$$
Let $Z+\xi\in (T\oplus T^*)\otimes\C$ and write
$$Z+\xi=Z^{1,0}+Z^{0,1}+\xi^{1,0}+\xi^{0,1}$$
for the decomposition with respect to the complex structure $I$ on $M$. A direct calculation shows that $Z+\xi$ lies in the $+i$-eigenspace $L_{\zeta}$ of $\mathcal{J}_{\zeta}$ if and only if
\begin{eqnarray*}
(1-|\zeta|^2)(-iZ^{1,0})+(-2\mathcal{I}m\zeta)(-\omega_J^{-1}\xi^{1,0})+(2\mathcal{R}e\zeta)(-\omega_K^{-1}\xi^{1,0}) & = & (1+|\zeta|^2)(+i)Z^{1,0}, \\
(1-|\zeta|^2)(+iZ^{0,1})+(-2\mathcal{I}m\zeta)(-\omega_J^{-1}\xi^{0,1})+(2\mathcal{R}e\zeta)(-\omega_K^{-1}\xi^{0,1}) & = & (1+|\zeta|^2)(+i)Z^{0,1}, \\
(1-|\zeta|^2)(+i\xi^{1,0})+(-2\mathcal{I}m\zeta)(\omega_JZ^{1,0})+(2\mathcal{R}e\zeta)(\omega_KZ^{1,0}) & = & (1+|\zeta|^2)(+i)\xi^{1,0}, \\
(1-|\zeta|^2)(-i\xi^{0,1})+(-2\mathcal{I}m\zeta)(\omega_JZ^{0,1})+(2\mathcal{R}e\zeta)(\omega_KZ^{0,1}) & = & (1+|\zeta|^2)(+i)\xi^{0,1}.
\end{eqnarray*}
The second and third equations are equivalent to the fourth and first equations respectively, which simplify to
\begin{eqnarray*}
\xi^{0,1} & = & i\mathcal{I}m\zeta(\omega_JZ^{0,1})-i\mathcal{R}e\zeta(\omega_KZ^{0,1}), \\
Z^{1,0} & = & -i\mathcal{I}m\zeta(\omega_J^{-1}\xi^{1,0})+i\mathcal{R}e\zeta(\omega_K^{-1}\xi^{1,0}).
\end{eqnarray*}
Substituting $\omega_J=(\sigma+\bar{\sigma})/2$ and $\omega_K=(\sigma-\bar{\sigma})/2i$, where $\sigma=t\sigma_Y$, and using the fact that $\sigma Z^{0,1}=0$ because $\sigma$ is of type $(2,0)$, yields
$$\xi^{0,1}=\frac{\zeta}{2}\bar{\sigma}Z^{0,1}$$
and similarly
$$Z^{1,0}=-\frac{\zeta}{2}\sigma^{-1}\xi^{1,0}.$$
This shows that $L_{\zeta}$ is the graph of the map
$$T^{0,1}\oplus \Omega^{1,0}\rightarrow \Omega^{0,1}\oplus T^{1,0}$$
given by the element
$$\frac{\zeta}{2}(-\sigma^{-1}+\bar{\sigma})=\frac{\zeta}{2}\left(-\frac{1}{t}\sigma_Y^{-1}+t\bar{\sigma}_Y\right)\in\Lambda^2T^{1,0}\oplus\Omega^{0,2}\subset \Lambda^2(\Omega^{0,1}\oplus T^{1,0}).$$
The result now follows from the theory of deformations of generalized complex structures (see Gualtieri~\cite{gualtieri11}).
\end{prf}

\begin{rmk}
We can also observe that the pure spinor defining the generalized complex structure $\mathcal{J}_{\zeta}$ is given by
$$\sigma_Y+2\zeta\left(\frac{1}{t}-\frac{1}{4}t\sigma_Y\bar{\sigma}_Y\right)-\zeta^2\bar{\sigma}_Y.$$
Applying the inverse of the isomorphism of Lemma~\ref{isomorphism} to the first order part in $\zeta$ yields $\frac{1}{2}\left(\frac{1}{t}\sigma_Y^{-1}-t\bar{\sigma}_Y\right)$. The direction of the deformation is then given by minus this class, as in Lemma~\ref{direction1}.
\end{rmk}

We finish this section with some discussion of Poisson deformations and $B$-field transforms of generalized complex manifolds. Suppose we are given a first order deformation, represented by
$$u_1\in C^{\infty}\left(\Lambda^2T^{1,0}\oplus (\Omega^{0,1}\otimes T^{1,0})\oplus\Omega^{0,2}\right).$$
To extend it to higher orders we must find a power series
$$u(\zeta)=u_1\zeta+u_2\zeta^2+u_3\zeta^3+\ldots$$
which satisfies the Maurer-Cartan equation
$$\bar{\partial}u(\zeta)+\frac{1}{2}[u(\zeta),u(\zeta)]=0,$$
where $[-,-]$ is the Schouten bracket. If
$$u_1\in C^{\infty}(\Lambda^2T^{1,0})$$
then the Maurer-Cartan equation decouples: the linear term and $C^{\infty}(\Lambda^3T^{1,0})$ part of the quadratic term give
$$\bar{\partial}u_1=0\qquad\mbox{and}\qquad [u_1,u_1]=0,$$
respectively. The first equation means that $u_1$ is a holomorphic section of $\Lambda^2T^{1,0}$ while the second means that $u_1$ is a Poisson structure. If $u_1$ satisfies these equations no higher order terms are needed, and convergence of $u(\zeta)=u_1\zeta$ is clearly automatic; we call this a holomorphic Poisson deformation (see Section~5.3 of Gualtieri~\cite{gualtieri11}).

For a K3 surface, $\sigma_Y^{-1}$ is a Poisson structure, with $[\sigma_Y^{-1},\sigma_Y^{-1}]=0$ following from $d\sigma_Y=0$. Moreover, $\sigma_Y^{-1}$ is everywhere of rank two. The Poisson deformation will therefore produce a generalized complex manifold which is everywhere of type $0$, meaning that it is a $B$-field transform of a symplectic manifold. 

The Maurer-Cartan equation also decouples when
$$u_1\in C^{\infty}(\Omega^{0,2}).$$
Moreover, we always have $[u_1,u_1]=0$, and for a K3 surface $\bar{\partial}u_1=0$ is automatic for degree reasons. Thus $u(\zeta)=u_1\zeta$ trivially satisfies the Maurer-Cartan equation and again no higher order terms are needed. This is known as a $B$-field transform, and we can explicitly write down the resulting $B$-field transform of the generalized complex structure $\mathcal{J}$.

By Lemma~\ref{direction2} the deformation $(M,\mathcal{J}_{\zeta})$ of $Y$ is in the direction
$$\frac{1}{2}\left(-\frac{1}{t}\sigma_Y^{-1}+t\bar{\sigma}_Y\right).$$
Choosing the right combination of Poisson deformation and $B$-field transform produces a family $(M,\mathcal{J}_{\zeta})$ that takes $Y$ to symplectic K3 surfaces, which is rather special behaviour. The moduli space of generalized K3 surfaces has complex dimension $22$. The $B$-field transforms of symplectic K3 surfaces form a dense subset, but the moduli space of genuine symplectic K3 surfaces has real dimension $22$; therefore we would not expect an arbitrary one-parameter family to contain any symplectic K3 surfaces at all.

\subsection{Deformations of categories and FM transforms}

We begin this section by reviewing Toda's work~\cite{toda09}. Let $X$ be a smooth projective variety over $\mathbb{C}$, and let $\mathrm{Coh}(X)$ be the category of coherent sheaves on $X$. First order deformations of $\mathrm{Coh}(X)$ as a $\mathbb{C}$-linear abelian category are parametrized by the degree two Hochschild cohomology $HH^2(X)$, which as we saw can be identified with $HT^2(X)$ using $I^{HKR}$. Toda described these deformations explicitly: given an element
$$u\in HT^2(X):=\H^0(X,\wedge^2T)\oplus \H^1(X,T)\oplus\H^2(X,\O_X)$$
he constructed a $\mathbb{C}[\epsilon]/(\epsilon^2)$-linear abelian category $\mathrm{Coh}(X,u)$. If $u$ lies in $\H^1(X,T)$ then $\mathrm{Coh}(X,u)$ arises from deforming $X$ as a complex manifold; if $u$ lies in $\H^0(X,\wedge^2T)$ then we get a ``non-commutative'' deformation; and if $u$ lies in $\H^2(X,\O_X)$ then we are led to a ``gerby'' deformation, consisting of twisted sheaves.

Toda also considered the behaviour of a Fourier-Mukai equivalence $\Phi:D^b(X)\rightarrow D^b(Y)$ under these deformations. As we have seen, $\Phi$ induces an isomorphism of Hochschild cohomology groups $\phi:HH^2(X)\rightarrow HH^2(Y)$. By composing this with the Hochschild-Kostant-Rosenberg isomorphism $I^{HKR}$, Toda {\em defined} an isomorphism $\phi_T:HT^2(X)\rightarrow HT^2(Y)$ (just as in Section~3.2). In this way, any element $u\in HT^2(X)$ corresponds to an element $v:=\phi_T(u)\in HT^2(Y)$, and we can consider the corresponding deformations $\mathrm{Coh}(X,u)$ and $\mathrm{Coh}(Y,v)$, or rather their derived categories $D^b(X,u)$ and $D^b(Y,v)$. Toda's main result is the existence of an equivalence
$$\Phi^{\dagger}:D^b(X,u)\rightarrow D^b(Y,v)$$
extending $\Phi$, in the sense that there is a commutative diagram
$$\xymatrix@R=20pt@C=20pt{
	D^b(X)\ar[rr]^{i_*}\ar[dd]^{\Phi} & & D^b(X,u)\ar[rr]^{Li^*}\ar[dd]^{\Phi^{\dagger}} & & D^-(X)\ar[dd]^{\Phi^-} \\
	 & & \\
	D^b(Y)\ar[rr]^{i_*} & & D^b(Y,v)\ar[rr]^{Li^*} & & D^-(Y). \\
	}$$
One interesting aspect is that $\phi_T$ need not preserve the summands of $HT^2$.

\begin{exm}
Let $X$ and $Y$ be dual complex tori, and $\Phi$ the Fourier-Mukai transform given by the Poincar{\'e} line bundle. Then $\phi_T$ maps $\H^1(X,T)$ to $\H^1(Y,T)$, because when we deform $X$ as a complex manifold it remains a complex torus, whose dual is the corresponding deformation of $Y$. On the other hand, Toda showed that $\phi_T$ maps $\H^0(X,\Lambda^2T)$ to $\H^2(Y,\O_Y)$ and maps $\H^2(X,\O_X)$ to $\H^0(Y,\Lambda^2T)$, meaning that non-commutative deformations of $\mathrm{Coh}(X)$ correspond to gerby deformations of $\mathrm{Coh}(Y)$, and vice versa. The resulting equivalence $\Phi^{\dagger}$ was extended to infinite order deformations by Ben-Bassat, Block, and Pantev~\cite{bbp07}.
\end{exm}

\begin{prf}{\bf of Theorem~\ref{main}(3)}
By Lemma~\ref{direction1} the family $(M,\mathcal{I}_{\zeta})$ is a deformation of $X$ in the direction
$$u_t=-\frac{2}{t}\sigma_X^{-1}[C_X]-\frac{2(t^2+1)}{t}\sigma_X^{-1}[F_X],$$
while by Lemma~\ref{direction2} the family $(M,\mathcal{J}_{\zeta})$ is a deformation of $Y$ in the direction
$$v_t=\frac{1}{2}\left(-\frac{1}{t}\sigma_Y^{-1}+t\bar{\sigma}_Y\right).$$
The Fourier-Mukai transform $\Phi: D^b(X)\rightarrow D^b(Y)$ induces an isomorphism $\phi_T:HT^2(X)\rightarrow HT^2(Y)$, which by Lemma~\ref{phi_T} takes $u_t$ to
\begin{eqnarray*}
\phi_T(u_t) & = & -\frac{2}{t}\left(\frac{1}{4}\sigma^{-1}_Y+\frac{2}{4}\bar{\sigma}_Y\right)-\frac{2(t^2+1)}{t}\left(-\frac{1}{4}\bar{\sigma}_Y\right) \\
 & = & \frac{1}{2}\left(-\frac{1}{t}\sigma_Y^{-1}+\frac{t^2-1}{t}\bar{\sigma}_Y\right) \\
 & = & v_t-\frac{1}{2t}\bar{\sigma}_Y.
\end{eqnarray*}
We have proved that $u_t$ corresponds to $v_t$ up to the $B$-field correction
$$-\frac{1}{2t}\bar{\sigma}_Y\in\H^2(Y,\O).$$
Moreover, the correction vanishes in the limit $t\rightarrow\infty$, establishing Theorem~\ref{main}(3).
\end{prf}

\begin{rmk}
If we use $\phi_{HT}$ instead of $\phi_T$, then by Lemma~\ref{phi'_T} we get an exact correspondence between $u_t$ and $v_t$, i.e., $v_t=\phi_{HT}(u_t)$ for all $t\in (1,\infty)$. Equivalently, $\phi_{H\Omega}$ takes the class
$$[\sigma_{\zeta}]=[\sigma_X]+2\zeta\left(\frac{1}{t}[C_X]+\frac{t^2+1}{t}[F_X]\right)-\zeta^2[\bar{\sigma}_X]$$
of the pure spinor of $(M,\mathcal{I}_{\zeta})$ to the class
$$[\sigma_Y]+2\zeta\left(\frac{1}{t}1_Y-\frac{1}{4}t[\sigma_Y\bar{\sigma}_Y]\right)-\zeta^2[\bar{\sigma}_Y]$$
of the pure spinor of $(M,\mathcal{J}_{\zeta})$. We don't know what this signifies; it is possible that $\phi_{HT}$ is the correct correspondence to use when considering deformations of generalized complex manifolds, as opposed to deformations of categories. This issue does not arise in the earlier example of dual complex tori, since the Todd class is trivial and $\phi_T=\phi_{HT}$ in that case.
\end{rmk}

\begin{cor}
\label{toda}
By Toda's results we obtain an equivalence of triangulated categories
$$\Phi^{\dagger}:D^b(X,u_t)\rightarrow D^b\left(Y,v_t-\frac{1}{2t}\bar{\sigma}_Y\right)$$
which extends $\Phi:D^b(X)\rightarrow D^b(Y)$.
\end{cor}

\begin{rmk}
An immediate question is whether the Fourier-Mukai transform $\Phi$ extends to higher orders, to infinite order (i.e., to categories over a formal neighbourhood of the point $0$ in $\P^1$ as in~\cite{bbp07}), or globally (i.e., to categories over $\P^1$). We will return to this question in Section~6.2.
\end{rmk}

\subsection{Limiting behaviour}

We are mainly interested in the limit $t\rightarrow\infty$, but let us first consider the limit $t\rightarrow 1^+$.

\begin{prp}
When $t=1$, the class $\alpha=[C_X]+2[F_X]$ on $X$ is semi-K{\"a}hler, by which we mean positive semi-definite. Nonetheless, we can still construct a twistor family of complex K3 surfaces $X_{\zeta}$, whose corresponding holomorphic two-forms are given by
$$\sigma_{\zeta}=\sigma_X+2\zeta\omega_I-\zeta^2\bar{\sigma}_X,$$
where the semi-K{\"a}hler form $\omega_I$ represents the class $\alpha$. As before, this is a deformation of $X$ in the direction
$$u_1=-2\sigma_X^{-1}[C_X]-4\sigma_X^{-1}[F_X].$$
\end{prp}

\begin{prf}
This class $\alpha$ is the limit of K{\"a}hler classes for $t>1$, but $\alpha.C_X=0$, so it must be on the wall of the K{\"a}hler cone. If we blow down the $(-2)$-curve $C_X$, we get an orbifold $\tilde{X}$ with a single $A_1$ singularity, admitting a K{\"a}hler class $\tilde{\alpha}$. We can apply the orbifold version of Yau's Theorem, which was proved by Kobayashi and Todorov~\cite{kt87}, building on unpublished ideas of Yau. (Note that the term ``generalized K3 surface'' means something quite different in~\cite{kt87}: it is used to denote a compact complex surface with at worst simple singular points whose minimal resolutions is a K3 surface.) This gives a Ricci-flat K{\"a}hler orbifold metric $\tilde{g}$ on $\tilde{X}$ whose K{\"a}hler form $\tilde{\omega}_I$ represents $\tilde{\alpha}$. The standard argument shows that $\tilde{g}$ is a hyperk{\"a}hler orbifold metric. Pulling back by $X\rightarrow\tilde{X}$ gives the required (singular) metric $g$ on $X$. The twistor family is then constructed in the usual way. 
\end{prf}

\begin{rmk}
For generic $\zeta$, $X_{\zeta}$ will be a non-algebraic K3 surface. Nevertheless, it comes from blowing up an orbifold, so it will always contain a complex curve, namely the exceptional divisor $C_X$.
\end{rmk}

\noindent
The corresponding deformation of $Y$ is in the direction
$$\phi_T(u_1)=-\frac{1}{2}\sigma_Y^{-1},$$
which gives a Poisson deformation. As discussed earlier, this deformation takes $Y$ to $B$-field transforms of symplectic K3 surfaces. It does not take $Y$ to a genuine symplectic K3 surface.

Next consider the limit $t\rightarrow\infty$. First note that the volume of a fibre of $X\rightarrow\P^1$ is $\alpha.F=\frac{1}{t}$, so $t\rightarrow\infty$ is the large complex structure limit where the fibres collapse (see Gross and Wilson~\cite{gw00}). We have
\begin{eqnarray*}
\lim_{t\rightarrow\infty}u_t & = & \lim_{t\rightarrow\infty}-\frac{2}{t}\sigma_X^{-1}[C_X]-\frac{2(t^2+1)}{t}\sigma_X^{-1}[F_X] \\
 & = & \lim_{t\rightarrow\infty}-2t\sigma_X^{-1}[F_X]
 \end{eqnarray*}
and
\begin{eqnarray*}
\lim_{t\rightarrow\infty}\phi_T(u_t) & = & \lim_{t\rightarrow\infty}v_t-\frac{1}{2t}\bar{\sigma}_Y \\
 & = & \lim_{t\rightarrow\infty}v_t \\
 & = & \lim_{t\rightarrow\infty}\frac{1}{2}\left(-\frac{1}{t}\sigma_Y^{-1}+t\bar{\sigma}_Y\right) \\
 & = & \lim_{t\rightarrow\infty}\frac{t}{2}\bar{\sigma}_Y.
\end{eqnarray*}
After renormalizing, it is clear that we should consider the deformation $X$ in the direction $-2\sigma_X^{-1}[F_X]$ and the corresponding deformation of $Y$ in the direction $\frac{1}{2}\bar{\sigma}_Y$. In fact, in these directions we obtain global deformations of not only the generalized complex manifolds, but also of the categories and Fourier-Mukai transform, as we now explain.

\begin{thm}
\label{gerby}
The corresponding directions
$$u_{\infty}:=-2\sigma_X^{-1}[F_X]\in\H^1(X,T)\qquad\mbox{and}\qquad v_{\infty}:=\phi_T(u_{\infty})=\frac{1}{2}\bar{\sigma}_Y\in\H^2(Y,\O)$$
produce one-parameter families of deformations of $X$ and $Y$, respectively, as generalized complex manifolds. The Fourier-Mukai transform $\Phi:D^b(X)\rightarrow D^b(Y)$ extends to a Fourier-Mukai transform between corresponding members of these families.
\end{thm}

\begin{prf}
We begin with $Y$. As explained earlier, when we deform $Y$ in the direction $v_{\infty}=\frac{1}{2}\bar{\sigma}_Y$, the resulting generalized K3 surface is a $B$-field transform of $Y$. When it comes to the derived category $D^b(Y)$, we ``integrate'' the first order deformation by exponentiating $v_{\infty}\zeta=\frac{1}{2}\bar{\sigma}_Y\zeta$ according to the map
$$\H^2(Y,\O)\stackrel{\exp}{\longrightarrow}\H^2(Y,\O^*)$$
to give a holomorphic gerbe $\beta_{\zeta}:=\exp(\frac{1}{2}\bar{\sigma}_Y\zeta)$ on $Y$. The resulting category is the derived category
$$D^b(Y,\beta_{\zeta}^{-1}):=D^b(\mathrm{Coh}(Y,\beta_{\zeta}^{-1}))$$
of $\beta_{\zeta}^{-1}$-twisted sheaves on $Y$ (as defined by C{\u a}ld{\u a}raru~\cite{caldararu00}). Note that Toda's convention for deforming $D^b(Y)$ includes a sign change for the $B$-field, which is why the category involves $\beta_{\zeta}^{-1}$ rather than $\beta_{\zeta}$.

The corresponding deformation of $X$ is a deformation as a complex manifold. By the Torelli theorem it suffices to describe the family of periods, given by the classes of the holomorphic two-forms. To first order in $\zeta$ these must look like $[\sigma_X]+2\zeta [F_X]$. In fact, no higher order terms are required since
$$([\sigma_X]+2\zeta [F_X])^2=[\sigma_X]^2+4\zeta [\sigma] [F_X]+ 4\zeta^2[F_X]^2=0.$$
We can also describe the deformations of $X$ more geometrically. First note that the class $[F_X]$ of a fibre is the pull-back of the generator $\eta_{\P^1}$ of $\H^2(\P^1,\mathbb{Z})$ under the map $\pi:X\rightarrow\P^1$. Thinking of $\eta_{\P^1}$ as a class in $\H^1(\P^1,\Omega^1_{\P^1})$, we can write
$$[F_X]=\pi^*\eta_{\P^1}\in\pi^*\H^1(\P^1,\Omega^1_{\P^1})\subset \H^1(X,\pi^*\Omega^1_{\P^1})\subset \H^1(X,\Omega^1_X).$$
The isomorphism $TX\cong\Omega^1_X$ induced by $\sigma_X$ takes $T_{X/\P^1}$ to $\pi^*\Omega^1_{\P^1}$, where $T_{X/\P^1}$ is the fibrewise tangent bundle, because $\pi:X\rightarrow\P^1$ is a Lagrangian fibration. Therefore
$$u_{\infty}=-2\sigma_X^{-1}[F_X]\in\H^1(X,T_{X/\P^1}).$$
This means that $u_{\infty}$ can be represented by a $1$-cocyle that on each $X_{ij}:=\pi^{-1}(U_i\cap U_j)$, where $\{U_i\}$ is an open cover of $\P^1$, is given by a vector field in the fibre direction. Moreover, these local vector field are constant in the fibre directions, since the class is pulled back from $\P^1$ (note that the open cover $\{U_i\}$ can be chosen so that no singular fibres lie inside the overlaps $X_{ij}$). Geometrically, this deformation produces a kind of torsor over $X$: integrating the  local vector fields produces translations in the fibre directions, and the resulting spaces $X_{\zeta}$ are obtained by taking the open sets $X_i:=\pi^{-1}(U_i)$ and gluing them on the overlaps $X_{ij}$ according to these translations. The deformation of $D^b(X)$ is of course given by $D^b(X_{\zeta})$.

Next we describe the equivalence of derived categories, which is an example of C{\u a}ld{\u a}raru's twisted Fourier-Mukai transforms~\cite{caldararu00}. Since $X_{\zeta}$ is a torsor over $X$, it has the same relative Jacobian (or dual fibration) as $X$, namely the elliptic surface $Y$. However, for general $\zeta$ the fibration $X_{\zeta}\rightarrow\P^1$ does not admit a section, and consequently there is no relative Poincar{\'e} sheaf on $X_{\zeta}\times_{\P^1} Y$. The existence of such a universal sheaf is obstructed by a holomorphic gerbe on $Y$, and this is precisely the gerbe $\beta_{\zeta}^{-1}$ from above. Instead, there is a $\pi_2^*\beta_{\zeta}^{-1}$-twisted universal sheaf on $X_{\zeta}\times_{\P^1} Y$, where $\pi_2:X_{\zeta}\times_{\P^1} Y\rightarrow Y$ is projection onto the second factor, and this leads to a twisted Fourier-Mukai transform
$$\Phi_{\zeta}:D^b(X_{\zeta})\rightarrow D^b(Y,\beta_{\zeta}^{-1}).$$

Finally, one can use the above equivalence to show that the one-parameter families of deformations of $X$ and $Y$ are in corresponding directions for all $\zeta$, not just at $\zeta=0$. The derived equivalence induces an isomorphism
$$\phi^{\zeta}_T:HT^2(X_{\zeta})\rightarrow HT^2(Y,\beta_{\zeta}^{-1}),$$
where on the right we consider a $\beta_{\zeta}^{-1}$-twisted version of $HT^2(Y)$ (cf.\ Huybrechts and Stellari~\cite{hs05}). Because $X_{\zeta}$ is only changing as a torsor over $X$, the deformation of $X_{\zeta}$ is in the direction
$$-2\sigma^{-1}_{X_{\zeta}}[F_{\zeta}]\in\H^1(X_{\zeta},T)\subset HT^2(X_{\zeta})$$
for all $\zeta$, where we use $[F_{\zeta}]$ to denote the class of a fibre of $X_{\zeta}\rightarrow\P^1$. On $Y$, only the gerbe changes, so the deformation of $(Y,\beta_{\zeta}^{-1})$ is in the direction
$$\frac{1}{2}\bar{\sigma}_Y\in\H^2(Y,\O)\cong \H^2(Y,\O,\beta_{\zeta}^{-1})\subset HT^2(Y,\beta_{\zeta}^{-1})$$
for all $\zeta$. One can check that these elements, $-2\sigma^{-1}_{X_{\zeta}}[F_{\zeta}]$ and $\frac{1}{2}\bar{\sigma}_Y$, correspond under the isomorphism $\phi^{\zeta}_T$.
\end{prf}

\begin{rmk}
The one-parameter families of the above theorem are families parametrized by $\C$, not by $\P^1$. In the case of $X$, the class of the holomorphic two-form $[\sigma_X]+2\zeta[F_X]$ tends to $[F_X]$ (up to scale) as $\zeta\rightarrow\infty$. Since $[F_x][\bar{F}_X]=[F_X]^2=0$ is not positive, $[F_X]$ cannot be the period of a complex K3 surface (really, we are approaching the boundary of the moduli space). Similarly, we know of no way to interpret the limit as $\zeta\rightarrow\infty$ of the gerby deformation of $Y$.
\end{rmk}

\begin{rmk}
The (analytic) Brauer group $\H^2(Y,\O^*)$ parametrizes torsors over the dual fibration $X$, as we can reconstruct $X_{\zeta}$ from $X\cong Y$ and $\beta_{\zeta}^{-1}\in\H^2(Y,\O^*)$. To understand the structure of the Brauer group, consider the long exact sequence
$$\ldots\rightarrow\H^2(Y,\Z)\stackrel{\iota}{\longrightarrow}\H^2(Y,\O)\stackrel{\exp}{\longrightarrow}\H^2(Y,\O^*)\rightarrow\H^3(Y,\Z)\rightarrow\ldots$$
coming from the exponential sequence. For K3 surfaces $\H^3(Y,\Z)$ is trivial and $\H^2(Y,\O)\cong\C$; therefore
$$\H^2(Y,\O^*)\cong \H^2(Y,\O)/\mathrm{image}(\iota)$$
is connected and one-dimensional. A generic elliptic K3 surface which admits a section will have Picard number $\rho=2$, and therefore the image of $\iota$ will have rank $22-\rho=20$. This means $\H^2(Y,\O^*)$ will be non-Hausdorff.
\end{rmk}

\section{Mirror symmetry}

\subsection{Differential-geometric}

Gross~\cite{gross99} established a differential-geometric version of mirror symmetry for K3 surfaces; we will follow the description in Section~1 of Gross and Wilson~\cite{gw00} (cf.\ also Proposition~6.8 of Huybrechts~\cite{huybrechts04}, which is based on work of Aspinwall and Morrison~\cite{am97}). Fix a sublattice of $\H^2(X,\mathbb{Z})$ isomorphic to the hyperbolic plane, with generators $[F]$ and $[C]$ satisfying
$$[F]^2=0,\qquad [C]^2=-2,\qquad\mbox{and}\qquad [F].[C]=1$$
(in our case these are the classes of the fibre and section of $X\rightarrow\P^1$, respectively). Then mirror symmetry is an involution acting on triples $(X,[B+i\omega],[\sigma])$, where $X$ is a marked K3 surface, $[B+i\omega]$ is a complexified K{\"a}hler class, and $[\sigma]$ is the class of a holomorphic two-form on $X$. The K{\"a}hler class $[\omega]$ should vanish on $F$, the class $[B]$ of the $B$-field should lie in $[F]^{\perp}/[F]\otimes\mathbb{R}$, and $\sigma$ should be normalized so that $[\mathcal{I}m\sigma]$ vanishes on $F$ and
$$[\omega]^2=[\mathcal{R}e\sigma]^2=[\mathcal{I}m\sigma]^2.$$
The mirror triple $({\check X},[{\check B}+i{\check \omega}],[{\check \sigma}])$ then satisfies
\begin{eqnarray*}
{[{\check \sigma}]} & \equiv & ([F].[\mathcal{R}e\sigma])^{-1}([C]+[B+i\omega])\qquad\mathrm{mod}[F], \\
{[{\check B}]} & \equiv & ([F].[\mathcal{R}e\sigma])^{-1}[\mathcal{R}e\sigma]-[C]\qquad\mathrm{mod}[F], \\
{[{\check \omega}]} & \equiv & ([F].[\mathcal{R}e\sigma])^{-1}[\mathcal{I}m\sigma]\qquad\mathrm{mod}[F],
\end{eqnarray*}
and is uniquely determined by the normalizations
$$[\mathcal{R}e\check \sigma]^2=[\mathcal{I}m\check \sigma]^2=[\check \omega]^2$$
and
$$[\check \omega].[\mathcal{R}e\check \sigma]=[\check \omega].[\mathcal{I}m\check \sigma]=[\mathcal{R}e\check \sigma].[\mathcal{I}m\check \sigma]=0.$$
Note that both $X$ and $\check X$ are complex K3 surfaces equipped with complexified K{\"a}hler classes, but one can simply ignore some of the data, namely $[B+i\omega]$ and $[\check\sigma]$, to obtain the more traditional mirror symmetry relation between a complex K3 surface $X$ and its mirror symplectic K3 surface $\check X$ (with $B$-field). This is the version we will use. Note that the above equations can then be summarized by
$$[\check B+i\check\omega]\equiv ([F].[\mathcal{R}e\sigma])^{-1}[\sigma] -[C]\qquad\mathrm{mod}[F].$$

\begin{prf}{\bf of Theorem~\ref{main}(4)}
Assume that $\zeta\neq 0$ or $\infty$. Recall that $(M,\mathcal{I}_{\zeta})$ is a complex K3 surface, with holomorphic two-form $\sigma_{\zeta}$ given by
$$\sigma_X+2\zeta\omega_I-\zeta^2\bar{\sigma}_X$$
up to scale, where $\omega_I$ represents the class $\alpha=\frac{1}{t}[C]+\left(\frac{t^2+1}{t}\right)[F]$. Since the normalization requires $[\mathcal{I}m\sigma_{\zeta}]$ to vanish on $F$, or equivalently, that $[\sigma_{\zeta}].[F]$ be real, we define
\begin{eqnarray*}
\sigma_{\zeta} & := & \frac{1}{2\zeta}(\sigma_X+2\zeta\omega_I-\zeta^2\bar{\sigma}_X) \\
 & = & \frac{\sigma_X}{2\zeta}+\omega_I-\frac{\zeta\bar{\sigma}_X}{2}.
\end{eqnarray*}
Then
$$[\sigma_{\zeta}].[F]=[\omega_I].[F]=\alpha.[F]=\frac{1}{t}$$
and according to the above relations, the mirror $(\check M,[\check B+i\check\omega],-)$ to $(M,-,[\sigma_{\zeta}])$ should satisfy
\begin{eqnarray*}
{[\check B+i\check\omega]} & \equiv & ([F].[\mathcal{R}e\sigma_{\zeta}])^{-1}[\sigma_{\zeta}] -[C]\qquad\mathrm{mod}[F] \\
 & \equiv & t\left(\frac{[\sigma_X]}{2\zeta}+\frac{1}{t}[C]+\left(\frac{t^2+1}{t}\right)[F]-\frac{\zeta[\bar{\sigma}_X]}{2}\right)-[C]\qquad\mathrm{mod}[F] \\
 & \equiv & \frac{t[\sigma_X]}{2\zeta}-\frac{\zeta t[\bar{\sigma}_X]}{2}\qquad\mathrm{mod}[F].
\end{eqnarray*}
By the remark at the end of Section~4, this is exactly the class of the complexified K{\"a}hler form on the symplectic K3 surface with $B$-field $(M,\mathcal{J}_{\zeta})$. We have shown that $(M,\mathcal{I}_{\zeta})=(M,-,[\sigma_{\zeta}])$ and $(M,\mathcal{J}_{\zeta})=(\check M,[\check B+i\check\omega],-)$ are mirror K3 surfaces for $\zeta\neq 0$ and $\infty$, in the sense of Gross. This establishes Theorem~\ref{main}(4), and completes the proof of the theorem.
\end{prf}

\subsection{Homological Mirror Symmetry}

We have shown that the complex K3 surface $(M,\mathcal{I}_{\zeta})$ and the symplectic K3 surface with $B$-field $(M,\mathcal{J}_{\zeta})$ are mirror partners in the sense of Gross~\cite{gross99} for $\zeta\neq 0$ and $\infty$. The Homological Mirror Symmetry (HMS) conjecture of Kontsevich~\cite{kontsevich95} asserts that the derived category of coherent sheaves on a complex manifold should be equivalent to the derived Fukaya category of the mirror symplectic manifold. In principle HMS is the strongest version of mirror symmetry, though it is not clear how to deduce from it other versions of mirror symmetry. HMS has been proved for quartic K3 surfaces by Seidel~\cite{seidel03}, but is currently an open conjecture for other K3 surfaces such as elliptic K3s. Theorem~\ref{main} suggests exploring HMS by trying to extend the Fourier-Mukai transform
$$\Phi:D^b(X)\rightarrow D^b(Y)$$
to some kind of equivalence for the families of generalized K3 surfaces $(M,\mathcal{I}_{\zeta})$ and $(M,\mathcal{J}_{\zeta})$. Since $X_{\zeta}=(M,\mathcal{I}_{\zeta})$ is a family of complex K3 surfaces, it is fairly clear that the appropriate deformation of $D^b(X)$ should just be the derived category $D^b(X_{\zeta})$ of coherent sheaves on $X_{\zeta}$.

\begin{que}
What is the appropriate deformation of $D^b(Y)$ associated to the family of generalized K3 surfaces $(M,\mathcal{J}_{\zeta})$? Does the derived category $D^b(Y)$ of coherent sheaves on the complex K3 surface $Y=(M,\mathcal{J}_0)$ deform to the derived Fukaya category of the symplectic K3 surface with $B$-field $(M,\mathcal{J}_{\zeta})$ when $\zeta\neq 0$?
\end{que}

For generic $\zeta$, the derived category $D^b(X_{\zeta})$ is relatively simple. In some sense, the complexity of a derived category can be measured by the number of spherical objects it contains. In~\cite{hsm08}, Huybrechts et al.\ studied derived categories of K3 surfaces and twisted K3 surfaces which contain at most one spherical object (up to shifts). A generic K3 surface $X_{\zeta}$ in a twistor family will be non-projective, and in fact will contain no complex curves, and so $\O_{X_{\zeta}}$ and its shifts will be the only spherical objects in $D^b(X_{\zeta})$. The limiting case $t=1$ is an exception, as then $X_{\zeta}$ always contains the $(-2)$-curve $C_X$; but this is a twistor family associated to a singular metric.

For the mirror manifold $(M,\mathcal{J}_{\zeta})$, the symplectic form
$$\omega=\csc\theta((\cos\phi)\omega_J+(\sin\phi)\omega_K)$$
and the $B$-field
$$B=\cot\theta((\sin\phi)\omega_J-(\cos\phi)\omega_K)$$
both vanish on the curve $C_Y\cong S^2$, which therefore yields a spherical object in the Fukaya category of $(M,\mathcal{J}_{\zeta})$. HMS should take $\O_{X_{\zeta}}$ to $C_Y$, or more generally, to a line bundle supported on $C_Y$ up to a shift (cf.\ the proof of Lemma~\ref{phiOmega}). Already we can start to surmise some correspondence between spherical objects.

Huybrechts, Marc{\`i}, and Stellari~\cite{hms09, hms11} also extended Toda's first order deformation theory for derived categories of K3 surfaces to infinite order, i.e., formal deformations.

\begin{que}
In Corollary~\ref{toda} we described a first order deformation
$$\Phi^{\dagger}:D^b(X,u_t)\rightarrow D^b\left(Y,v_t-\frac{1}{2t}\bar{\sigma}_Y\right)$$
of the Fourier-Mukai transform
$$\Phi:D^b(X)\rightarrow D^b(Y).$$
What is the significance of the $B$-field correction $\frac{1}{2t}\bar{\sigma}_Y$? Does $\Phi^{\dagger}$ extend to infinite order?
\end{que}

More generally, one could try to associate to an arbitrary generalized complex manifold some analogue of the derived category of coherent sheaves, and then develop generalized Fourier-Mukai transforms. Presumably the category should be built from generalized complex branes, which are described by Gualtieri~\cite{gualtieri11}. These include holomorphic bundles and Lagrangian submanifolds with flat bundles when the generalized complex manifold is of complex or symplectic type, respectively, and so there is some hope of building a category from generalized complex branes that would specialize to the derived category of coherent sheaves or to the derived Fukaya category in these two cases.

\begin{flushleft}
Department of Mathematics\hfill sawon@email.unc.edu\\
University of North Carolina\hfill www.unc.edu/$\sim$sawon\\
Chapel Hill NC 27599-3250\\
USA\\
\end{flushleft}


\begin{thebibliography}{XXX}



\bibitem{am97} P.\ Aspinwall and D.\ Morrison,
{\em String theory on K3 surfaces},
Mirror symmetry II, 703--716, AMS/IP Stud.\ Adv.\ Math.\ {\bf 1}, Amer.\ Math.\ Soc., Providence, RI, 1997.

\bibitem{bbhm98} C.\ Bartocci, U.\ Bruzzo, D.\ Hern{\'a}ndez Ruip{\'e}rez, and J.M.\ Mu{\~n}oz Porras,
{\em Mirror symmetry on K3 surfaces via Fourier-Mukai transform},
Comm.\ Math.\ Phys.\ {\bf 195} (1998), no.\ 1, 79--93.





\bibitem{bbp07} O.\ Ben-Bassat, J.\ Block, and T.\ Pantev,
{\em Non-commutative tori and Fourier-Mukai duality},
Compos.\ Math.\ {\bf 143} (2007), no.\ 2, 423--475.


\bibitem{crv09} D.\ Calaque, C.\ Rossi, and M.\ Van den Bergh,
{\em C{\u a}ld{\u a}raru's conjecture and Tsygan's formality},
preprint {\bf arXiv:0904.4890}.

\bibitem{caldararu00} A.\ C{\u a}ld{\u a}raru,
{\em Derived categories of twisted sheaves on Calabi-Yau manifolds},
Cornell PhD thesis, May 2000 (available from {\bf
 www.math.upenn.edu/$\sim$andreic/}).


\bibitem{caldararu03i} A.\ C{\u a}ld{\u a}raru,
{\em The Mukai pairing, I: the Hochschild structure},
preprint {\bf math.AG/0308079}.

\bibitem{caldararu05} A.\ C{\u a}ld{\u a}raru,
{\em The Mukai pairing, II: the Hochschild-Kostant-Rosenberg isomorphism},
Adv.\ Math.\ {\bf 194} (2005), no.\ 1, 34--66.

\bibitem{cw10} A.\ C{\u a}ld{\u a}raru and S.\ Willerton,
{\em The Mukai pairing, I: a categorical approach},
New York J.\ Math.\ {\bf 16} (2010), 61--98.
















\bibitem{goto05} R.\ Goto,
{\em On deformations of generalized Calabi-Yau, hyperK{\"a}hler, $G_2$ and Spin$(7)$ structures I},
preprint {\bf arXiv:math/0512211}.



\bibitem{gross99} M.\ Gross,
{\em Special Lagrangian fibrations, II: geometry. A survey of techniques in the study of special Lagrangian fibrations},
in Surveys in differential geometry: differential geometry inspired by string theory, 341--403,
Surv.\ Differ.\ Geom.\ {\bf 5}, Int.\ Press, Boston, MA, 1999.

\bibitem{gw00} M.\ Gross and P.\ Wilson,
{\em Large complex structure limits of K3 surfaces},
J.\ Differential Geom.\ {\bf 55} (2000), no.\ 3, 475--546.


\bibitem{gualtieri11} M.\ Gualtieri,
{\em Generalized complex geometry},
Ann.\ of Math.\ (2) {\bf 174} (2011), no.\ 1, 75--123.







\bibitem{hitchin03} N.\ Hitchin,
{\em Generalized Calabi-Yau manifolds},
Q.\ J.\ Math.\ {\bf 54} (2003), no.\ 3, 281--308.

\bibitem{hklr87} N.\ Hitchin, A.\ Karlhede, U.\ Lindstr{\"o}m, and M.\ Ro{\v c}ek,
{\em Hyper-K{\"a}hler metrics and supersymmetry},
Comm.\ Math.\ Phys.\ {\bf 108} (1987), no.\ 4, 535--589.









\bibitem{huybrechts04} D.\ Huybrechts,
{\em Moduli spaces of hyperk{\"a}hler manifolds and mirror symmetry},
Intersection theory and moduli, 185--247 (electronic),
ICTP Lect.\ Notes {\bf XIX}, Abdus Salam Int.\ Cent.\ Theoret.\ Phys., Trieste, 2004.




\bibitem{hsm08} D.\ Huybrechts, E.\ Macr{\`i}, and P.\ Stellari,
{\em Stability conditions for generic K3 categories},
Compos.\ Math.\ {\bf 144} (2008), no.\ 1, 134--162.

\bibitem{hms09} D.\ Huybrechts, E.\ Macr{\`i}, and P.\ Stellari,
{\em Derived equivalences of K3 surfaces and orientation},
Duke Math.\ J.\ {\bf 149} (2009), no.\ 3, 461--507.

\bibitem{hms11} D.\ Huybrechts, E.\ Macr{\`i}, and P.\ Stellari,
{\em Formal deformations and their categorical general fibre},
Comment.\ Math.\ Helv.\ {\bf 86} (2011), no.\ 1, 41--71.

\bibitem{hn11} D.\ Huybrechts and M.\ Nieper-Wi{\ss}kirchen,
{\em Remarks on derived equivalences of Ricci-flat manifolds},
Math.\ Z.\ {\bf 267} (2011), no.\ 3--4, 939--963.

\bibitem{hs05} D.\ Huybrechts and P.\ Stellari,
{\em Equivalences of twisted K3 surfaces},
Math.\ Ann.\ {\bf 332} (2005), no.\ 4, 901--936.





\bibitem{kt87} R.\ Kobayashi and A.\ Todorov,
{\em Polarized period map for generalized K3 surfaces and the moduli of Einstein metrics},
Tohoku Math.\ J.\ (2) {\bf 39} (1987), no.\ 3, 341--363.

\bibitem{kontsevich95} M.\ Kontsevich,
{\em Homological algebra of mirror symmetry},
Proceedings of the International Congress of Mathematicians, Vol.\ 1, 2 (Z{\"u}rich, 1994), 120--139, Birkh{\"a}user, Basel, 1995.



\bibitem{mns08} E.\ Macr{\`i}, M.\ Nieper-Wi{\ss}kirchen, and P.\ Stellari,
{\em The module structure of Hochschild homology in some examples},
C.\ R.\ Math.\ Acad.\ Sci.\ Paris {\bf 346} (2008), no.\ 15--16, 863--866.

\bibitem{ms09} E.\ Macr{\`i} and P.\ Stellari,
{\em Infinitesimal derived Torelli theorem for K3 surfaces}
(Appendix by S.\ Mehrotra),
Int.\ Math.\ Res.\ Not.\ IMRN 2009, 3190--3220.











\bibitem{mukai87} S.\ Mukai,
{\em Fourier functor and its application to the moduli of bundles on an abelian variety},
Advanced Studies in Pure Mathematics {\bf 10} (1987), Algebraic Geometry, Sendai, 1985, 515--550.










\bibitem{ramadoss08} A.\ Ramadoss,
{\em The relative Riemann-Roch theorem from Hochschild homology},
New York J.\ Math.\ {\bf 14} (2008), 643--717.










\bibitem{seidel03} P.\ Seidel,
{\em Homological mirror symmetry for the quartic surface},
arXiv preprint {\bf math/0310414}.


\bibitem{toda09} Y.\ Toda,
{\em Deformations and Fourier-Mukai transforms},
J.\ Differential Geom.\ {\bf 81} (2009), no.\ 1, 197--224.











\end{thebibliography}
\end{document}